\documentclass{amsart}

\usepackage{amsfonts,amssymb,amscd,amsmath}

\usepackage[dvips]{graphicx}

\usepackage{psfrag,euscript}

\usepackage[latin1]{inputenc}

\newtheorem{theorem}{Theorem}
\theoremstyle{plain}

\newtheorem{conjecture}{Conjecture}
\newtheorem{corollary}{Corollary}

\newtheorem{lemma}{Lemma}

\newtheorem{proposition}{Proposition}

\numberwithin{equation}{section}

\newcommand{\ov}{\overline}

\newcommand{\leb}{\operatorname{Leb}}

\newcommand{\dist}{\operatorname{dist}}

\def \CC {{\mathbb C}}

\def \KK {{\mathbb K}}
\def \LL {{\mathbb L}}

\def \PP {{\mathbb P}}

\def \RR {{\mathbb R}}

\def \ZZ {{\mathbb Z}}

\newcommand{\cF}{\EuScript{F}}

\newcommand{\cD}{\EuScript{D}}

\newcommand{\cO}{\EuScript{O}}

\def \cD {{\mathcal D}}

\def \cF {{\mathcal F}}

\def \cK {{\mathcal K}}

\def \cO {{\mathcal O}}

\def \cU {{\mathcal U}}
\def \cV {{\mathcal V}}

\def \fX {{\mathfrak X}}

\begin{document}

\title[Sensitive dependence vs. physical measure]
{On sensitive dependence on initial conditions and
  existence of physical measure for $3$-flows}

\author{Vítor Araújo}
\address{V\'{\i}tor Ara\'ujo, Instituto de Matem\'a\-tica,
Universidade Federal da Bahia\\
Av. Adhemar de Barros, S/N , Ondina,
40170-110 - Salvador-BA-Brazil}

\email{vitor.d.araujo@ufba.br}

\begin{abstract}After reviewing known results on sensitiveness and
  also on robustness of attractors together with
  observations on their proofs, we show that for attractors
  of three-dimensional flows, robust chaotic behavior
  (meaning sensitiveness to initial conditions for the past
  as well for the future for all nearby flows) is equivalent
  to the existence of certain hyperbolic structures. These
  structures, in turn, are associated to the existence of
  physical measures. In short \emph{in low dimensions,
    robust chaotic behavior for smooth flows ensures the
    existence of a physical measure}.
\end{abstract}

\thanks{
Author was partially supported by FAPERJ (Rio de
Janeiro), FAPESB (Bahia), CNPq and
PRONEX-Dynamical Systems (Brazil).
}

\subjclass{Primary: 37D25; Secondary: 37D30, 37D45.}
\renewcommand{\subjclassname}{\textup{2000} Mathematics Subject Classification}
\keywords{sensitive dependence on initial conditions,
  physical measure, singular-hyperbolicity, expansiveness,
  robust attractor, robust chaotic flow, positive Lyapunov
  exponent}

\date{2012}

\maketitle
\tableofcontents

\section{Introduction}
\label{sec:introd}

The development of the theory of dynamical systems has shown
that models involving expressions as simple as quadratic
polynomials (as the \emph{logistic family} or \emph{H\'enon
  attractor}, see e.g.\cite{devaney1989} for a gentle
introduction), or autonomous ordinary differential equations
with a hyperbolic equilibrium of saddle-type accumulated by
regular orbits, as the \emph{Lorenz flow} (see
e.g. \cite{GH83,viana2000i,AraPac2010}), exhibit \emph{sensitive
  dependence on initial conditions}, a common feature of
\emph{chaotic dynamics}: small initial differences are
rapidly augmented as time passes, causing two trajectories
originally coming from practically indistinguishable points
to behave in a completely different manner after a short
while. Long term predictions based on such models are
unfeasible since it is not possible to both specify initial
conditions with arbitrary accuracy and numerically calculate
with arbitrary precision. For an introduction to these
notions see~\cite{devaney1989,robinson2004}.

\begin{figure}[htpb]
\psfrag{X}{$X^t(x)$}\psfrag{Y}{$X^t(y)$}
\psfrag{x}{$x$}\psfrag{y}{$y$}
\includegraphics[width=4cm]{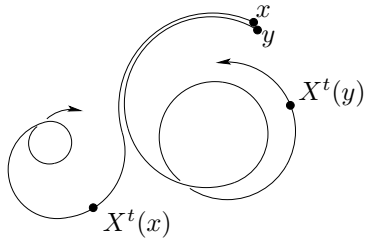}
\caption{\label{fig-sensivel} Sensitive dependence on
  initial conditions.}
\end{figure}

Formally the definition of sensitivity is as follows for a
flow $X^t$ on some compact manifold $M$: a $X^t$-invariant
subset $\Lambda$ is \emph{sensitive to initial conditions}
or has \emph{sensitive dependence on initial conditions}, or
simply \emph{chaotic} if, for every small enough $r>0$ and
$x\in\Lambda$, and for any neighborhood $U$ of $x$, there
exists $y\in U$ and $t\neq0$ such that $X^t(y)$ and $X^t(x)$
are $r$-apart from each other:
$\dist\big(X^t(y),X^t(x)\big)\ge r$. See
Figure~\ref{fig-sensivel}. An analogous definition holds for
diffeomorphism $f$ of some manifold, taking $t\in\ZZ$ and
setting $f=X^1$ in the previous definition.

Using some known results on robustness of attractors from
Ma\~n\'e~\cite{Man82} and Morales, Pacifico and
Pujals~\cite{MPP04} together with observations on their
proofs, we show that for attractors of three-dimensional
flows, robust chaotic behavior (in the above sense of
sensitiveness to initial conditions) is equivalent to the
existence of certain hyperbolic structures. These
structures, in turn, are associated to the existence of
physical measures. In short \emph{in low dimensions, robust
  chaotic behavior ensures the existence of a physical
  measure}.

\section{Preliminary notions}
\label{sec:prelim-notions}

Here and throughout the text we assume that $M$ is a
three-dimensional compact connected manifold without
boundary endowed with some Riemannian metric which induces a
distance denoted by $\dist$ and a volume form $\leb$ which
we name \emph{Lebesgue measure} or \emph{volume}. For any
subset $A$ of $M$ we denote by $\ov{A}$ the (topological)
closure of $A$.

We denote by $\fX^r(M), r\ge1$ the set of $C^r$ smooth
vector fields $X$ on $M$ endowed with the $C^r$
topology. Given $X\in\fX^r(M)$ we denote by $X^t$, with
$t\in\RR$, the flow generated by the vector field $X$. Since
we assume that $M$ is a compact manifold the flow is defined
for all time. Recall that the flow
$(X^t)_{t\in\RR}$ is a family of $C^r$ diffeomorphisms
satisfying the following properties:
\begin{enumerate}
\item $X^0=Id:M\to M$ is the identity map of $M$;
\item $X^{t+s}=X^t\circ X^s$ for all $t,s\in\RR$,
\end{enumerate}
and it is \emph{generated by the vector field} $X$ if
\begin{enumerate}
\item[3] $\left.\frac{d}{dt} X^t (q)\right|_{t=t_0} =
X\big(X_{t_0}(q)\big)$ for all $q\in M$ and $t_0\in\RR$.
\end{enumerate}

We say that a compact $X^t$-invariant set $\Lambda$ is
\emph{isolated} if there exists a neighborhood $U$ of
$\Lambda$ such that $\Lambda=\cap_{t\in\RR}X^t(U)$.  A
compact invariant set $\Lambda$ is {\em attracting} if
$\Lambda_X(U):=\cap_{t\geq 0}X^t(U)$ equals $\Lambda$ for
some neighborhood $U$ of $\Lambda$ satisfying
$\ov{X^t(U)}\subset U$, for all $t> 0$. In this case the
neighborhood $U$ is called an \emph{isolating} neighborhood
of $\Lambda$. Note that $\Lambda_X(U)$ is in general
different from $\cap_{t\in\RR} X^t(U)$, but for an
attracting set the extra condition $\ov{X^t(U)}\subset U$
for $t>0$ ensures that every attracting set is also
isolated. We say that $\Lambda$ is {\em transitive} if
$\Lambda$ is the closure of both $\{X^t(q):t>0\}$ and
$\{X^t(q):t<0\}$ for some $q\in \Lambda$.  An {\em
  attractor} of $X$ is a transitive attracting set of $X$
and a {\em repeller} is an attractor for $-X$.  We say that
$\Lambda$ is a \emph{proper} attractor or repeller if
$\emptyset\neq\Lambda\neq M$.

An \emph{equilibrium} (or \emph{singularity}) for $X$ is a
point $\sigma\in M$ such that $X^t(\sigma)=\sigma$ for all
$t\in\RR$, i.e.  a fixed point of all the flow maps, which
corresponds to a zero of the associated vector field $X$:
$X(\sigma)=0$. An \emph{orbit} of $X$ is a set
$\cO(q)=\cO_X(q)=\{X^t(q): t\in\RR\}$ for some $q\in M$.  A
{\em periodic orbit} of $X$ is an orbit $\cO=\cO_X(p)$ such
that $X^T(p)=p$ for some minimal $T>0$. A \emph{critical
  element} of a given vector field $X$ is either an
equilibrium or a periodic orbit.

We recall that a $X^t$-invariant probability measure $\mu$
is a probability measure satisfying $\mu(X^t(A))=\mu(A)$ for
all $t\in\RR$ and measurable $A\subset M$. Given an
invariant probability measure $\mu$ for a flow $X^t$, let
$B(\mu)$ be the the \emph{(ergodic) basin} of $\mu$, i.e.,
the set of points $z\in M$ satisfying for all continuous
functions $\varphi: M \to \RR$
\begin{align*}
  \lim_{T\to+\infty} \frac{1}{T} \int_0^T \varphi\big(
  X^t(z) \big) \, dt = \int\varphi\,d\mu.
\end{align*}
We say that $\mu$ is a \emph{physical} (or \emph{SRB})
measure for $X$ if $B(\mu)$ has positive Lebesgue measure:
$\leb\big( B(\mu) \big)>0$. 

The existence of a physical measures for an attractor shows
that most points in a neighborhood of the attractor have
well defined long term statistical behavior. So, in spite of
chaotic behavior preventing the exact prediction of the time
evolution of the system in practical terms, we gain some
statistical knowledge of the long term behavior of the
system near the chaotic attractor.

\section{Chaotic systems}
\label{sec:chaotic-systems}

We distinguish between forward and backward sensitive
dependence on initial conditions. We say that an invariant
subset $\Lambda$ for a flow $X^t$ is \emph{future chaotic
  with constant $r>0$} if, for every $x\in\Lambda$ and each
neighborhood $U$ of $x$ in the ambient manifold, there
exists $y\in U$ and $t>0$ such that
$\dist\big(X^t(y),X^t(x)\big)\ge r$. Analogously we say that
$\Lambda$ is \emph{past chaotic with constant $r$} if
$\Lambda$ is future chaotic with constant $r$ for the flow
generated by $-X$. If we have such \emph{sensitive
  dependence both for the past and for the future}, we say
that $\Lambda$ is \emph{chaotic}. Note that in this language
sensitive dependence on initial conditions is weaker than
chaotic, future chaotic or past chaotic conditions.

An easy consequence of chaotic behavior is that it prevents
the existence of sources or sinks, either attracting or
repelling equilibria or periodic orbits, inside the
invariant set $\Lambda$. Indeed, if $\Lambda$ is future
chaotic (for some constant $r>0$) then, were it to contain
some attracting periodic orbit or equilibrium, any point of
such orbit (or equilibrium) would admit no point in a
neighborhood whose orbit would move away in the
future. Likewise, reversing the time direction, a past
chaotic invariant set cannot contain repelling periodic
orbits or repelling equilibria. As an almost reciprocal we
have the following.

\begin{lemma}
  \label{le:nonchaotic-interio}
  If $\Lambda=\cap_{t\in\RR}X^t(U)$ is a compact isolated
  proper subset for $X\in\fX^1(M)$ with isolating
  neighborhood $U$ and $\Lambda$ is \emph{not future
    chaotic} (respective not past chaotic), then
  $\Lambda^-_X(U):=\overline{\cap_{t>0}X^{-t}(U)}$
  (respective $\Lambda^+_X(U):=\overline{\cap_{t>0}X^t(U)}$)
  has non-empty interior.
\end{lemma}

\begin{proof}
  If $\Lambda$ is not future chaotic, then for every $r>0$
  there exists some point $x\in\Lambda$ and a neighborhood
  $V$ of $x$ such that $\dist\big( X^t(y), X^t(x)\big)<r$
  for all $t>0$ and each $y\in V$. If we choose
  $0<r<\dist(M\setminus U,\Lambda)$ (we note that if
  $\Lambda=U$ then $\Lambda$ would be open and closed, and
  so, by connectedness of $M$, $\Lambda$ would not be a
  proper subset), then we deduce that $X^t(y)\in U$, that
  is, $y\in X^{-t}(U)$ for all $t>0$, hence $V\subset
  \Lambda^-_X(U)$. Analogously if $\Lambda$ is not past
  chaotic, just by reversing the time direction.  
\end{proof}

In particular \emph{if an invariant and isolated set $\Lambda$
  with isolating neighborhood $U$ is given such that the volume of
  both $\Lambda_X^+(U)$ and $\Lambda^-_X(U)$ is zero, then
  $\Lambda$ is chaotic}.

Sensitive dependence on initial conditions is part of the
definition of \emph{chaotic system} in the literature, see
e.g.~\cite{devaney1989}. It is an interesting fact that
sensitive dependence is a consequence of another two common
features of most systems considered to be chaotic: existence
of a dense orbit and existence of a dense subset of periodic
orbits. 

\begin{proposition}
  \label{pr:transdensperiodic-sensitive}
  A compact invariant subset $\Lambda$ for a flow $X^t$ with
  a dense subset of periodic orbits and a dense (regular and
  non-periodic) orbit is chaotic.
\end{proposition}

A short proof of this proposition can be found
in~\cite{BBCDS92}. An extensive discussion of this and
related topics can be found in~\cite{GlasWeiss93}.

\section{Lack of sensitiveness for flows on surfaces}
\label{sec:flows-surfac}

We recall the following celebrated result of Mauricio
Peixoto in \cite{peixoto59,peixoto62} (and for a more
detailed exposition of this results and sketch of the proof
see \cite{GH83}) built on previous work of Poincar\'e
\cite{PoincareI} and Andronov and
Pontryagin \cite{AP37}, that characterizes structurally
stable vector fields on compact surfaces.
\begin{theorem}[Peixoto]
  \label{thm:Peixoto}
  A $C^r$ vector field, $r\ge1$, on a compact surface $S$ is
  structurally stable if, and only if:
  \begin{enumerate}
  \item the number of critical elements is
    finite and each is hyperbolic;
  \item there are no orbits connecting saddle points;
  \item the non-wandering set consists of critical elements
    alone.
  \end{enumerate}
  Moreover if $S$ is orientable, then the set of
  structurally stable vector fields is open and dense in
  $\fX^r(S)$.
\end{theorem}
In particular, this implies that for a structurally stable
vector field $X$ on $S$ there is an open and dense subset
$B$ of $S$ such that the positive orbit $X^t(p), t\ge0$ of $p\in
B$ converges to one of finitely many attracting equilibria.
Therefore \emph{no sensitive dependence on initial conditions
arises for an open and dense subset of all vector fields in
orientable surfaces}.

The extension of Peixoto's characterization of structural
stability for $C^r$ flows, $r\ge1$, on non-orientable
surfaces is known as \emph{Peixoto's Conjecture}, and up
until now it has been proved for the projective plane $\PP^2$
\cite{PM82}, the Klein bottle $\KK^2$ \cite{Markley69} and
$\LL^2$, the torus with one cross-cap \cite{gutierrez78}.
Hence for these surfaces we also have no sensitiveness to
initial conditions for most vector fields.

This explains in part the great interest attached to the
Lorenz attractor which was one of the first examples of
sensitive dependence on initial conditions.

\section{Robustness and volume hyperbolicity}
\label{sec:robustn-singul-hyper}

Related to chaotic behavior is the notion of \emph{robust
  dynamics}. We say that an attracting set
$\Lambda=\Lambda_X(U)$ for a $3$-flow $X$ and some open
subset $U$ is \emph{robust} if there exists a $C^1$
neighborhood $\cU$ of $X$ in $\fX^1(M)$ such that
$\Lambda_Y(U)$ is transitive for every $Y\in\cU$. 

The following result obtained by Morales, Pacifico and
Pujals in~\cite{MPP04} characterizes robust attractors for
three-dimensional flows.
\begin{theorem}
  \label{thm:robust-attractor-sing-hyp}
  Robust attractors for flows containing equilibria are
  singular-hyperbolic sets for $X$.
\end{theorem}
We remark that robust attractors cannot be $C^1$
approximated by vector fields presenting either attracting
or repelling periodic points.  This implies that, on
$3$-manifolds, any periodic orbit inside a robust
attractor is hyperbolic of saddle-type.

We now define the concept of singular-hyperbolicity.
A compact invariant set $\Lambda$ of $X$ is {\em partially
  hyperbolic} if there are a continuous invariant tangent
bundle decomposition $T_\Lambda M=E^s_\Lambda\oplus
E^c_\Lambda$ and constants $\lambda,K>0$ such that
\begin{itemize}
\item
{\em $E^c_\Lambda$ $(K,\lambda)$-dominates $E^s_\Lambda$},
i.e. for all $x\in\Lambda$ and for all $t\ge0$
\begin{align}\label{eq.domination}
  \| DX^t(x)\mid E^s_x\| \leq \frac{e^{-\lambda t}}{K}\cdot
  m(DX^t(x) \mid E^c_x);
\end{align}
\item
$E^s_\Lambda$ is $(K,\lambda)$-contracting: $\|DX^t\mid
E^s_x\|\le K e^{-\lambda t}$ for all $x\in\Lambda$ and for all $t\ge0$.
\end{itemize}
For $x\in \Lambda$ and $t\in\RR$ we let $J_t^c(x)$ be
the absolute value of the determinant of the linear map
$DX^t(x)\mid E^c_x:E^c_x\to E^c_{X^t(x)}$.  We say that the
sub-bundle $E^c_\Lambda$ of the partial hyperbolic set
$\Lambda$ is {\em $(K,\lambda)$-volume expanding} if
$$
J_t^c(x)=\big| \det(DX^t\mid E^c_x) \big|\geq Ke^{\lambda t},
$$
for every $x\in \Lambda$ and $t\geq 0$.

We say that a partially hyperbolic set is {\em
  singular-hyperbolic} if its singularities are hyperbolic
and it has volume expanding central direction.

A \emph{singular-hyperbolic attractor} is a
singular-hyperbolic set which is an attractor as well: an
example is the (geometric) Lorenz attractor presented in
\cite{Lo63,Gu76}. Any equilibrium $\sigma$ of a
singular-hyperbolic attractor for a vector field $X$ is such
that $DX(\sigma)$ has only real eigenvalues
$\lambda_2\le\lambda_3\le\lambda_1$ satisfying the same
relations as in the Lorenz flow example:
\begin{align}
  \label{eq:lorenz-like}
  \lambda_2<\lambda_3<0<-\lambda_3<\lambda_1,
\end{align}
which we refer to as \emph{Lorenz-like equilibria}.
We recall that an compact $X^t$-invariant set $\Lambda$ is
\emph{hyperbolic} if the tangent bundle over $\Lambda$
splits $T_\Lambda M=E^s_\Lambda\oplus E^X_\Lambda \oplus
E^u_\Lambda$ into three $DX^t$-invariant subbundles, where
$E^s_\Lambda$ is uniformly contracted, $E^u_\Lambda$ is
uniformly expanded, and $E^X_\Lambda$ is the direction of
the flow at the points of $\Lambda$.  It is known,
see~\cite{MPP04}, that a partially hyperbolic set for a
three-dimensional flow, with volume expanding central
direction and without equilibria, is hyperbolic. Hence the
notion of singular-hyperbolicity is an extension of the
notion of hyperbolicity.

Recently in a joint work with Pacifico, Pujals and Viana
\cite{APPV} the following consequence of transitivity and
singular-hyperbolicity was proved.

\begin{theorem}
  \label{srb}
  Let $\Lambda=\Lambda_X(U)$ be a singular-hyperbolic
  \emph{attractor} of a flow $X\in\fX^2(M)$ on a
  three-dimensional manifold.  Then $\Lambda$ supports a
  unique physical probability measure $\mu$ which is ergodic
  and its ergodic basin covers a full Lebesgue measure
  subset of the topological basin of attraction, i.e.
  $B(\mu)=W^s(\Lambda)$ Lebesgue mod $0$.  Moreover the
  support of $\mu$ is the whole attractor $\Lambda$.
\end{theorem}
It follows from the proof in~\cite{APPV} that the
singular-hyperbolic \emph{attracting set} $\Lambda_Y(U)$ for
all $Y\in\fX^2(M)$ which are $C^1$-close enough to $X$
\emph{admits finitely many physical measures whose ergodic
  basins cover $U$ except for a zero volume subset}.

\subsection{Absence of sinks and sources nearby}
\label{sec:absence-sinks-source}

The proof of Theorem~\ref{thm:robust-attractor-sing-hyp}
given in~\cite{MPP04} uses several tools from the theory of
normal hyperbolicity developed first by Ma\~n\'e
in~\cite{Man82} together with the low dimension of the
flow. Indeed, going through the proof in~\cite{MPP04} we can
see that the arguments can be carried through assuming that
\begin{enumerate}
\item $\Lambda$ is an attractor for $X$ with isolating
  neighborhood $U$ such that every equilibria in $U$ is
  hyperbolic with no resonances;
\item there exists a $C^1$ neighborhood $\cU$ of $X$ such
  that for all $Y\in\cU$ every periodic orbit and
  equilibria in $U$ is hyperbolic of saddle-type.
\end{enumerate}
The condition on the equilibria amounts to restricting
the possible three-dimensional vector fields in the above
statement to an open a dense subset of all $C^1$ vector
fields. Indeed, the hyperbolic and no-resonance condition on
a equilibrium $\sigma$ means that:
\begin{itemize}
\item either $\lambda\neq\Re(\omega)$ if the eigenvalues of
  $DX(\sigma)$ are $\lambda\in\RR$ and
  $\omega,\overline\omega\in\CC$;
\item or $\sigma$ has only real eigenvalues with different
  norms.
\end{itemize}
Indeed, conditions (1) and (2) ensure that no bifurcations
of periodic orbits or equilibria leading to sinks or sources
are allowed for any nearby flow in $U$. This implies, by now
standard arguments, that the flow on $\Lambda$ must have a
dominated splitting which is \emph{volume hyperbolic}: both
subbundles of the splitting must contract/expand volume. For
a $3$-dimensional flow one of the subbundles is
one-dimensional, and so we deduce singular-hyperbolicity
either for $X$ or for $-X$.  If $\Lambda$ has no equilibria,
then $\Lambda$ is uniformly hyperbolic. Otherwise, it
follows from the arguments in~\cite{MPP04} that all
singularities of $\Lambda$ are Lorenz-like and this shows
that $\Lambda$ must be singular-hyperbolic for $X$.

We note that the second condition above is a consequence of
any one of the following assumptions on $U$:
\begin{description}
\item[\textbf{robust chaoticity}] for every $Y\in\cU$ the maximal
  invariant subset $\Lambda_Y(U)$ is chaotic;
\item[\textbf{zero volume and future chaoticity}] for every
  $Y\in\cU$ the maximal invariant subset $\Lambda_Y(U)$ has
  zero volume and is future chaotic;
\item[\textbf{zero volume and robust positive Lyapunov
  exponent}] for every $Y\in\cU$ the maximal
  invariant subset $\Lambda_Y(U)$ has zero volume and
  there exists a full Lebesgue measure subset $P_Y$ of $U$
  such that
  \begin{align}
    \label{eq:positiveLyap}
    \limsup\frac1n\sum_{i=0}^{n-1}\log\|DY^i_x\|>0, \quad
    x\in P_Y.
  \end{align}
\end{description}
The result of  Ma\~n\'e analogous to
Theorem~\ref{thm:robust-attractor-sing-hyp} in~\cite{Man82}
\begin{theorem}
  \label{thm:robust2d-hyp}
  Robust attractors for surface diffeomorphisms are hyperbolic.
\end{theorem}
also follows from the absence of sinks and sources for all
$C^1$ close diffeomorphisms in a neighborhood of the
attractor.

Extensions of these results to higher dimensions for
diffeomorphisms, by Bonatti, Díaz and Pujals in~\cite{BDP},
show that robust transitive sets always admit a volume
hyperbolic splitting of the tangent bundle.  Vivier in
\cite{Viv03} extends previous results of Doering~\cite{Do87}
for flows, showing that a $C^1$ robustly transitive vector
field on a compact boundaryless $n$-manifold, with $n\ge3$,
admits a global dominated splitting. Metzger and Morales
extend the arguments in \cite{MPP04} to homogeneous vector
fields (inducing flows allowing no bifurcation of critical
elements, i.e. no modification of the index of periodic
orbits or equilibria) in higher dimensions leading to the
concept of $2$-sectional expanding attractor
in~\cite{MeMor06}.

\subsection{Robust chaoticity, volume hyperbolicity and
  physical measure}
\label{sec:robust-chaotic-parti}

The preceding observations allows us to deduce that robust
chaoticity is a sufficient conditions for
singular-hyperbolicity of a generic attractor.

\begin{corollary}
  \label{cor:robust-chaotic-sing-hyp}
  Let $\Lambda$ be an attractor for $X\in\fX^1(M^3)$ such
  that every equilibrium in its trapping region is
  hyperbolic with no resonances. Then $\Lambda$ is
  singular-hyperbolic if, and only if, $\Lambda$ is robustly
  chaotic.
\end{corollary}

This means that \emph{if we can show that arbitrarily close
  orbits, in an isolating neighborhood of an attractor, are
  driven apart, for the future as well as for the past, by
  the evolution of the system, and this behavior persists
  for all $C^1$ nearby vector fields, then the attractor is
  singular-hyperbolic}.

To prove the necessary condition on
Corollary~\ref{cor:robust-chaotic-sing-hyp} we use the
concept of expansiveness for flows, and through it show that
singular-hyperbolic attractors for $3$-flows are robustly
expansive and, as a consequence, robustly chaotic also. This
is done in the last Section~\ref{sec:expans-systems}.

We recall the following conjecture of Viana, presented in~\cite{Vi98}
\begin{conjecture}\label{conj:viana}
  If an attracting set $\Lambda(U)$ of smooth map/flow has a
  non-zero Lyapunov exponent at Lebesgue almost every point
  of its isolated neighborhood $U$ (i.e. it
  satisfies~\eqref{eq:positiveLyap} with $P_Y\subset U$),
  then it admits some physical measure.
\end{conjecture}
From the preceding results and observations we can give a
partial answer to this conjecture for $3$-flows in the
following form.

\begin{corollary}
\label{cor:weakconjecture}
Let $\Lambda_X(U)$ be an attractor for a flow $X\in\fX^1(M)$
such that
\begin{itemize}
\item the divergence of $X$ is negative in $U$;
\item the equilibria in $U$ are hyperbolic with no resonances;
\item there exists a neighborhood $\cU$ of $X$ in $\fX^1(M)$
  such that for $Y\in\cU$ one has \eqref{eq:positiveLyap}
  almost everywhere in $U$.
\end{itemize}
Then there exists a neighborhood $\cV\subset\cU$ of $X$ in
$\fX^1(M)$ and a dense subset $\cD\subset\cV$ such that
\begin{enumerate}
\item  $\Lambda_Y(U)$ is singular-hyperbolic for all $Y\in\cV$;
\item there exists a physical measure $\mu_Y$ supported in
  $\Lambda_Y(U)$ for all $Y\in\cD$.
\end{enumerate}
\end{corollary}
Indeed, item (2) above is a consequence of item (1), the
denseness of $\fX^2(M)$ in $\fX^1(M)$ in the $C^1$ topology,
together with Theorem~\ref{srb} and the observation
following its statement.

Item (1) above is a consequence of
Corollary~\ref{cor:robust-chaotic-sing-hyp} and the
observations of Section~\ref{sec:absence-sinks-source},
noting that negative divergence on the isolating
neighborhood $U$ ensures that the volume of $\Lambda_Y(U)$
is zero for $Y$ in a $C^1$ neighborhood $\cV$ of $X$.

\section{Expansive systems}
\label{sec:expans-systems}

Here we explain why robust chaotic behavior necessarily
follows from singular-hyperbolicity in an attractor,
completing the proof of
Corollary~\ref{cor:robust-chaotic-sing-hyp}. For this we
need the concept and some properties of expansiveness for
flows.

A concept related to sensitiveness is that of expansiveness,
which roughly means that points whose orbits are always
close for all time must coincide.  The concept of
expansiveness for homeomorphisms plays an important role in
the study of transformations. Bowen and
Walters~\cite{BoWa72} gave a definition of expansiveness for
flows which is now called \emph{C-expansiveness}
cite{KS79}. The basic idea of their definition is that two
points which are not close in the orbit topology induced by
$\RR$ can be separated at the same time even if one allows a
continuous time lag --- see below for the technical
definitions.  The equilibria of C-expansive flows must be
isolated \cite[Proposition 1]{BoWa72} which implies that the
Lorenz attractors and geometric Lorenz models are not
C-expansive.

Keynes and Sears introduced \cite{KS79} the idea of
restriction of the time lag and gave several definitions of
expansiveness weaker than C-expansiveness. The notion of
\emph{K-expansiveness} is defined allowing only the time lag
given by an increasing surjective homeomorphism of $\RR$.
Komuro \cite{Km84} showed that the Lorenz attractor and the
geometric Lorenz models are not K-expansive.  The reason for
this is not that the restriction of the time lag is
insufficient, but that the topology induced by $\RR$ is
unsuited to measure the closeness of two points in the same
orbit.

Taking this fact into consideration, Komuro \cite{Km84} gave
a definition of \emph{expansiveness} suitable for flows
presenting equilibria accumulated by regular orbits. This
concept is enough to show that two points which do not lie
on a same orbit can be separated.

Let $C\big(\RR,\RR\big)$ be the set of all continuous
functions $h:\RR\to\RR$ and let us write
$C\big((\RR,0),(\RR,0)\big)$ for the subset of all $h\in
C\big(\RR,\RR\big)$ such that $h(0)=0$. We define
\begin{align*}
  \cK_0=\{ h\in C\big((\RR,0),(\RR,0)\big): h(\RR)=\RR, \,
  h(s)>h(t)\,,\forall s>t \},
\end{align*}
and
\begin{align*}
  \cK=\{ h\in C\big(\RR,\RR\big) : h(\RR)=\RR, \,
  h(s)>h(t)\,,\forall s>t \},
\end{align*}

A flow $X$ is \emph{C-expansive} (\emph{K-expansive}
respectively) on an invariant subset $\Lambda\subset M$ if
for every $\epsilon>0$ there exists $\delta>0$ such that if
$x,y\in \Lambda$ and for some $h\in
\cK_0$ (respectively $h\in\cK$\,) we have
\begin{align}\label{eq:distorbit}
  \dist\big( X^t(x) , X^{h(t)}(y) \big) \le \delta
  \quad\text{for all}\quad t\in\RR,
\end{align}
then $y\in X^{[-\epsilon,\epsilon]}(x)=\{ X^t(x):
-\epsilon\le t \le\epsilon\}$.

We say that the flow $X$ is \emph{expansive} on $\Lambda$ if
for every $\epsilon>0$ there is $\delta >0$ such that for
$x,y\in \Lambda$ and some $h\in\cK$ (note that now we do not
demand that $0$ be fixed by $h$)
 satisfying
\eqref{eq:distorbit}, then we can find $t_0 \in \RR$ such
that $X^{h(t_0)}(y) \in X^{[t_0-\epsilon,t_0+\epsilon]}(x)$.

Observe that expansiveness on $M$ implies sensitive
dependence on initial conditions for any flow on a manifold
with dimension at least 2. Indeed if $\epsilon,\delta$
satisfy the expansiveness condition above with $h$ equal to
the identity and we are given a point $x\in M$ and a
neighborhood $U$ of $x$, then taking $y\in U\setminus
X^{[-\epsilon,\epsilon]}(x)$ (which always exists since we
assume that $M$ is not one-dimensional) there exists
$t\in\RR$ such that $\dist\big(X^t(y),X^t(x)\big)\ge\delta$.
The same argument applies whenever we have expansiveness on
an $X$-invariant subset $\Lambda$ of $M$ containing a dense
regular orbit of the flow.

Clearly C-expansive $\implies$ K-expansive $\implies$
expansive by definition. When a flow has no fixed point then
C-expansiveness is equivalent to K-expansiveness
\cite[Theorem A]{Oka90}. In \cite{BoWa72} it is shown that
on a connected manifold a C-expansive flow has no fixed
points. The following was kindly communicated to us by
Alfonso Artigue from IMERL, the proof can be found
in~\cite{AraPac2010}.

\begin{proposition}
  \label{pr:KCequiv}
  A flow is C-expansive on a manifold $M$ if, and only if,
  it is K-expansive.
\end{proposition}

We will see that singular-hyperbolic attractors are
expansive. In particular, the Lorenz attractor and the
geometric Lorenz examples are all expansive and sensitive to
initial conditions. Since these families of flows exhibit
equilibria accumulated by regular orbits, we see that
expansiveness is compatible with the existence of fixed
points by the flow.

\subsection{Singular-hyperbolicity and expansiveness}
\label{sec:singul-hyperb-expans}

The full proof of the following given in~\cite{APPV}. We
provide a sketch of the proof in
Section~\ref{sec:proof-expansiveness}.

\begin{theorem}
  \label{thm:sing-hyp-attract-expansive}
Let $\Lambda$ be a singular-hyperbolic attractor of $X\in
{\fX}^1(M)$.  Then $\Lambda$ is expansive.
\end{theorem}



The reasoning is based on analyzing Poincar\'e return maps
of the flow to a convenient ($\delta$-adapted)
cross-section.  We use the family of adapted cross-sections
and corresponding Poincará maps $R$, whose Poincar\'e time
$t(\cdot)$ is large enough, obtained assuming that the
attractor $\Lambda$ is singular-hyperbolic. These cross-sections
have a co-dimension $1$ foliation, which are dynamically
defined, whose leaves are uniformly contracted and invariant
under the Poincaré maps.  In addition $R$ is uniformly
expanding in the transverse direction and this also holds
near the singularities.

From here we argue by contradiction: if the flow is not
expansive on $\Lambda$, then we can find a pair of orbits
hitting the cross-sections infinitely often on pairs of
points uniformly close. We derive a contradiction by showing
that the uniform expansion in the transverse direction to
the stable foliation must take the pairs of points apart,
unless one orbit is on the stable manifold of the other.

This argument only depends on the existence of finitely many
Lorenz-like singularities on a compact partially hyperbolic
invariant attracting subset $\Lambda=\Lambda_X(U)$, with
volume expanding central direction, and of a family of
adapted cross-sections with Poincaré maps between them,
whose derivative is hyperbolic. It is straightforward that
if these conditions are satisfied for a flow $X^t$ of
$X\in\fX^1(M^3)$, then the same conditions hold for all
$C^1$ nearby flows $Y^t$ and for the maximal invariant
subset $\Lambda_Y(U)$ \emph{with the same family of
  cross-sections} which are also adapted to $\Lambda_Y(U)$
(as long as $Y$ is $C^1$-close enough to $X$).

\begin{corollary}
  \label{cor:robustexpansive}
  A singular-hyperbolic attractor $\Lambda=\Lambda_X(U)$ is
  \emph{robustly expansive}, that is, there exists a
  neighborhood $\cU$ of $X$ in $\fX^1(M)$ such that
  $\Lambda_Y(U)$ is expansive for each $Y\in\cU$, where $U$
  is an isolating neighborhood of $\Lambda$.
\end{corollary}

Indeed, since transversality, partial hyperbolicity and
volume expanding central direction are robust properties,
and also the hyperbolicity of the Poincaré maps depends on
the central volume expansion, all we need to do is to check
that a given adapted cross-section $\Sigma$ to $X$
is also adapted to $Y\in\fX^1$ for every $Y$ sufficiently
$C^1$ close to $X$. But $\Lambda_X(U)$ and $\Lambda_Y(U)$
are close in the Hausdorff distance if $X$ and $Y$ are close
in the $C^0$ distance, by the following elementary result.

\begin{lemma}
  \label{le:upper-semicont-maximal}
  Let $\Lambda$ be an isolated set of $X\in {\fX}^r(M)$,
  $r\geq 0$.  Then for every isolating block $U$ of
  $\Lambda$ and every $\epsilon>0$ there is a neighborhood
  ${\cU}$ of $X$ in ${\fX}^r(M)$ such that
  $\Lambda_Y(U)\subset B(\Lambda,\epsilon)$ and
  $\Lambda\subset B(\Lambda_Y(U),\epsilon)$ for all
  $Y\in{\cU}$.
\end{lemma}

Thus, if $\Sigma$ is an adapted cross-section we can find a
$C^1$-neighborhood $\cU$ of $X$ in $\fX^1$ such that
$\Sigma$ is still adapted to every flow $Y^t$ generated by a
vector field in $\cU$.

\subsection{Singular-hyperbolicity and chaotic behavior}
\label{sec:singul-hyperb-chaoti}

We already know that expansiveness implies sensitive
dependence on initial conditions. An argument with the same
flavor as the proof of expansiveness provides the following,
whose proof also sketch in the following
Section~\ref{sec:proof-expansiveness}. See also
\cite{AMS2010} for a different approach to sensitiveness.

\begin{theorem}\label{thm:sing-hyp-chaotic}
  A singular-hyperbolic isolated set $\Lambda=\cap_{t\in\RR}
  \overline{X^t(U)}$ is \emph{robustly chaotic}, i.e. there
  exists a neighborhood $\cU$ of $X$ in $\fX^1(M)$ such that
  $\cap_{t\in\RR}\overline{Y^t(U)}$ is chaotic for each
  $Y\in\cU$, where $U$ is an isolating neighborhood of
  $\Lambda$.
\end{theorem}

This completes the argument proving that robust chaoticity
is a necessary property of singular-hyperbolicity, in
Corollary~\ref{cor:robust-chaotic-sing-hyp}.

\section{Sketch of the proof of expansiveness and of chaotic
behavior}
\label{sec:proof-expansiveness}

We need the following notions to understand the proof of
Theorems~\ref{thm:sing-hyp-attract-expansive} and
Theorem~\ref{thm:sing-hyp-chaotic} as a consequence.


\subsection{Adapted cross-sections and Poincar\'e maps}
\label{sec:cross-sect-poinc}

To help explain the ideas of the proofs we give here a few
properties of \emph{Poincar\'e maps}, that is, continuous
maps $R:\Sigma\to\Sigma'$ of the form $R(x)=X^{t(x)}(x)$
between cross-sections $\Sigma$ and $\Sigma'$ of the flow
near a singular-hyperbolic set.  We always assume that the
Poincar\'e time $t(\cdot)$ is large enough as explained in
what follows.

We assume that $\Lambda$ is a compact invariant subset for a
flow $X\in\fX^1(M)$ such that $\Lambda$ is a
singular-hyperbolic attractor.  Then every equilibrium in
$\Lambda$ is Lorenz-like.


\subsubsection{Stable foliations on cross-sections}
\label{s.21}

We start recalling standard facts about uniformly hyperbolic
flows from e.g. \cite{HPS77}.

An embedded disk $\gamma\subset M$ is a (local) {\em
  strong-unstable manifold}, or a {\em strong-unstable
  disk}, if $\dist(X^{-t}(x),X^{-t}(y))$ tends to zero
exponentially fast as $t\to+\infty$, for every
$x,y\in\gamma$. Similarly, $\gamma$ is called a (local) {\em
  strong-stable manifold}, or a {\em strong-stable disk}, if
$\dist(X^{t}(x),X^{t}(y))\to0$ exponentially fast as
$n\to+\infty$, for every $x,y\in\gamma$. It is well-known
that every point in a uniformly hyperbolic set possesses a
local strong-stable manifold $W_{loc}^{ss}(x)$ and a local
strong-unstable manifold $W_{loc}^{uu}(x)$ which are disks
tangent to $E_x$ and $G_x$ at $x$ respectively with
topological dimensions $d_E=\dim(E)$ and $d_G=\dim(G)$
respectively. Considering the action of the flow we get the
(global) \emph{strong-stable manifold}
$$W^{ss}(x)=\bigcup_{t>0}
X^{-t}\Big(W^{ss}_{loc}\big(X^t(x)\big)\Big)$$
and the
(global) \emph{strong-unstable manifold}
$$W^{uu}(x)=\bigcup_{t>0}X^{t}\Big(W^{uu}_{loc}\big(X_{-t}(x)\big)\Big)$$
for every point $x$ of a uniformly hyperbolic set.  These
are immersed submanifolds with the same differentiability of
the flow.  We also consider the \emph{stable manifold}
$W^s(x)=\cup_{t\in\RR} X^{t}\big(W^{ss}(x)\big)$ and
\emph{unstable manifold}
$W^u(x)=\cup_{t\in\RR}X^{t}\big(W^{uu}(x)\big)$ for $x$ in
a uniformly hyperbolic set, which are flow invariant.

Now we recall classical facts about partially hyperbolic
systems, especially existence of strong-stable and
center-unstable foliations.  The standard reference is
\cite{HPS77}.

We have that $\Lambda$ is a singular-hyperbolic isolated set
of $X\in {\fX}^1(M)$ with invariant splitting $T_\Lambda M =
E^{s}\oplus E^{cu}$ with $\dim E^{cu}=2$.  Let
$\tilde{E}^s\oplus \tilde{E}^{cu}$ be a continuous extension
of this splitting to a small neighborhood $U_0$ of
$\Lambda$.  For convenience we take $U_0$ to be forward
invariant.  Then $\tilde{E}^s$ may be chosen invariant under
the derivative: just consider at each point the direction
formed by those vectors which are strongly contracted by
$DX^t$ for positive $t$.  In general $\tilde{E}^{cu}$ is not
invariant. However we can consider a cone field around it on
$U_0$
$$
C^{cu}_a(x)=\{v=v^s+v^{cu}: v^s\in \tilde{E}^s_x
\text{ and }
v^{cu}\in\tilde{E}^{cu}_x
\text{ with } \|v^s\|\le a\cdot \|v^{cu}\|\}
$$
which is forward invariant for $a>0$:
\begin{equation}
\label{eq.cone3}
DX^t(C^{cu}_a(x)) \subset C^{cu}_a (X^t(x))
\quad\text{for all large $t>0$.}
\end{equation}
Moreover we may take $a>0$ arbitrarily small, reducing
$U_0$ if necessary.  For notational simplicity we write
$E^s$ and $E^{cu}$ for $\tilde E^s$ and $\tilde E^{cu}$ in
all that follows.

From the standard normal hyperbolic theory, there are
locally strong-stable and center-unstable manifolds, defined
at every regular point $x\in U_0$ and which are embedded
disks tangent to $E^s(x)$ and $E^{cu}(x)$, respectively. 
 Given any
$x\in U_0$ define the strong-stable manifold $W^{ss}(x)$ and
the stable-manifold $W^s(x)$ as for an hyperbolic flow
(see the beginning of this section).

Denoting $E^{cs}_x=E^s_x\oplus E^{X}_x$, where $E^X_x$ is
the direction of the flow at $x$, it follows that
$$
T_x W^{ss}(x)=E^s_x \quad\text{and}\quad
T_x W^{s}(x)=E^{cs}_x.
$$
We fix $\epsilon$ once and for all. Then we call
$W^{ss}_{\epsilon}(x)$ the local \emph{strong-stable manifold}
and $W^{cu}_{\epsilon}(x)$ the local \emph{center-unstable
manifold} of $x$.

Now let $\Sigma$ be a \emph{cross-section} to the flow, that
is, a $C^1$ embedded compact disk transverse to $X$: at
every point $z\in\Sigma$ we have $T_z\Sigma\oplus E^X_z=T_z
M$ (recall that $E^X_z$ is the one-dimensional subspace
$\{s\cdot X(z): s\in\RR\}$).  For every $x\in\Sigma$ we
define $W^s(x,\Sigma)$ to be the connected component of
$W^s(x)\cap\Sigma$ that contains $x$.  This defines a
foliation $\cF^{s}_{\Sigma}$ of $\Sigma$ into co-dimension
$1$ sub-manifolds of class $C^1$.

\begin{figure}[ht]
\centering
\psfrag{a}{$W^s(x,\sigma)$}
\psfrag{b}{$R_\Sigma$}
\psfrag{c}{$W^s(x)$}
\psfrag{d}{$W^{ss}(x)$}
\psfrag{e}{$\Sigma(\gamma)$}
\psfrag{f}{$\Sigma$}
\psfrag{g}{$\gamma$}
  \includegraphics[width=7cm,height=4.5cm]{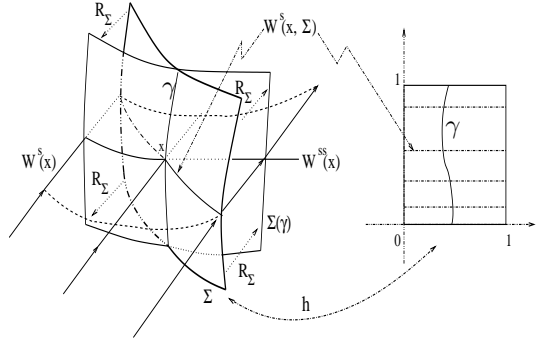}
  \caption{\label{f.squaresection} The sections $\Sigma$,
    $\Sigma(\gamma)$, the manifolds $W^s(x), W^{ss}(x)$,
    $W^s(x,\Sigma)$ and the projection $R_\Sigma$, on the
    right. On the left, the square $[0,1]\times[0,1]$ is
    identified with $\Sigma$ through the map $h$, where
    $\cF_\Sigma^s$ becomes the horizontal foliation and the
    curve $\gamma$ is transverse to the horizontal
    direction. Solid lines with arrows indicate the flow
    direction.}
\end{figure}

Given any cross-section $\Sigma$ and a point $x$ in its
interior, we may always find a smaller cross-section also
with $x$ in its interior and which is the image of the
square $[0,1]\times[0,1]$ by a $C^1$ diffeomorphism $h$ that
sends horizontal lines inside leaves of $\cF^{s}_{\Sigma}$.
In what follows we assume that cross-sections are of this
kind, see Figure~\ref{f.squaresection}.

We also assume that each cross-section $\Sigma$ is contained
in $U_0$, so that every $x\in\Sigma$ is such that
$\omega(x)\subset \Lambda$.

On the one hand $x\mapsto W^{ss}_\epsilon(x)$ is usually not
differentiable if we assume that $X$ is only of class $C^1$,
see e.g. \cite{PT93}.  On the other hand, assuming that the
cross-section is small with respect to $\epsilon$, and
choosing any curve $\gamma\subset\Sigma$ crossing
transversely every leaf of $\cF_\Sigma^s$\,, we may consider
a Poincar\'e map
$$
R_\Sigma:\Sigma \to \Sigma(\gamma)= \bigcup_{z\in\gamma}
W^{ss}_\epsilon(z)
$$
with Poincar\'e time close to zero, see
Figure~\ref{f.squaresection}.  This is a homeomorphism onto
its image, close to the identity, such that
$R_\Sigma(W^s(x,\Sigma))\subset
W^{ss}_\epsilon(R_\Sigma(x))$.  So, identifying the points
of $\Sigma$ with their images under this homeomorphism, we
pretend that indeed $W^s(x,\Sigma)\subset
W^{ss}_\epsilon(x)$.

\subsubsection{Hyperbolicity of Poincar\'e maps}
\label{s.22}

Let $\Sigma$ be a small cross-section to $X$ and let
$R:\Sigma\to\Sigma'$ be a Poincar\'e map $R(y)=X^{t(y)}(y)$ to
another cross-section $\Sigma'$ (possibly $\Sigma=\Sigma'$).
Here $R$ needs not correspond to the first time the orbits
of $\Sigma$ encounter $\Sigma'$.

The splitting $E^s\oplus E^{cu}$ over $U_0$ induces a continuous
splitting $E_\Sigma^s\oplus E_\Sigma^{cu}$ of the tangent bundle
$T\Sigma$ to $\Sigma$ (and analogously for $\Sigma'$), defined by
\begin{equation}\label{eq.splitting}
E_\Sigma^s(y)=E^{cs}_y\cap T_y{\Sigma}
\quad\mbox{and}\quad
E_\Sigma^{cu}(y)=E^{cu}_y\cap T_y{\Sigma}.
\end{equation}
We now show that if the Poincar\'e time $t(x)$ is
sufficiently large then \eqref{eq.splitting} defines a
hyperbolic splitting for the transformation $R$ on the
cross-sections restricted to $\Lambda$.

\begin{proposition}\label{p.secaohiperbolica}
Let $R:\Sigma\to\Sigma'$ be a Poincar\'e map as before
with Poincar\'e time $t(\cdot)$.
Then $DR_x(E_\Sigma^s(x)) = E_{\Sigma'}^s(R(x))$ at every
$x\in\Sigma$ and $DR_x(E_\Sigma^{cu}(x)) = E_{\Sigma'}^{cu}(R(x))$
at every $x\in\Lambda\cap\Sigma$.

Moreover for every given $0<\lambda<1$ there exists
$T_1=T_1(\Sigma,\Sigma',\lambda)>0$ such that if
$t(\cdot)>T_1$ at every point, then
$$
\|DR \mid E^s_\Sigma(x)\| < \lambda
\quad\text{and}\quad
\|DR \mid E^{cu}_\Sigma(x)\| > 1/\lambda
\quad\text{at every $x\in\Sigma$.}
$$
\end{proposition}

Given a cross-section $\Sigma$, a positive number $\rho$,
and a point $x\in \Sigma$,
we define the unstable cone of width $\rho$ at $x$ by
\begin{equation}
\label{cone}
C_\rho^u(x)=\{v=v^s+v^u : v^s\in E^s_\Sigma(x),\, v^u\in E^{cu}_\Sigma(x)
\mbox{ and } \|v^s\| \le \rho \|v^u\| \}.
\end{equation}

Let $\rho>0$ be any small constant. In the following
consequence of Proposition~\ref{p.secaohiperbolica} we
assume the neighborhood $U_0$ has been chose sufficiently
small.
\begin{corollary}
\label{ccone}
For any $R:\Sigma\to\Sigma'$ as in
Proposition~\ref{p.secaohiperbolica}, with $t(\cdot)>T_1$\,,
and any $x \in\Sigma$, we have $DR(x) (C^u_{\rho}(x))
\subset C_{\rho/2}^u(R(x))$ and
\begin{align*}
  \| DR_x(v)\| \ge \frac5{6}\lambda^{-1} \cdot \|v\|
  \quad\mbox{for all}\quad v\in C^u_{\rho}(x).
\end{align*}
\end{corollary}
As usual a \emph{curve} is the image of a compact interval
$[a,b]$ by a $C^1$ map. We use $\ell(\gamma)$ to denote its
length.  By a \emph{cu-curve} in $\Sigma$ we mean a curve
contained in the cross-section $\Sigma$ and whose tangent
direction is contained in the unstable cone
$T_z\gamma\subset C^u_\rho(z)$ for all $z\in\gamma$.  The
next lemma says that \emph{the length of cu-curves linking the stable
  leaves of nearby points $x,y$ must be bounded by the
  distance between $x$ and $y$}.
\begin{lemma}
\label{l.lengthversusdistance}
Let us we assume that $\rho$ has been fixed, sufficiently small.
Then there exists a constant $\kappa$ such that, for any pair
of points $x, y \in \Sigma$, and any cu-curve $\gamma$ joining $x$ to
some point of $W^s(y,\Sigma)$, we have $\ell(\gamma)\le\kappa
\cdot d(x,y)$.
\end{lemma}
Here $d$ is the intrinsic distance in the $C^2$ surface
$\Sigma$, that is, the length of the shortest smooth curve
inside $\Sigma$ connecting two given points in $\Sigma$.

In what follows we take $T_1$ in
Proposition~\ref{p.secaohiperbolica} for $\lambda=1/3$.


\subsubsection{Adapted cross-sections}
\label{s.23}
\begin{figure}[htpb]
\centering
    \includegraphics[height=3.5cm]{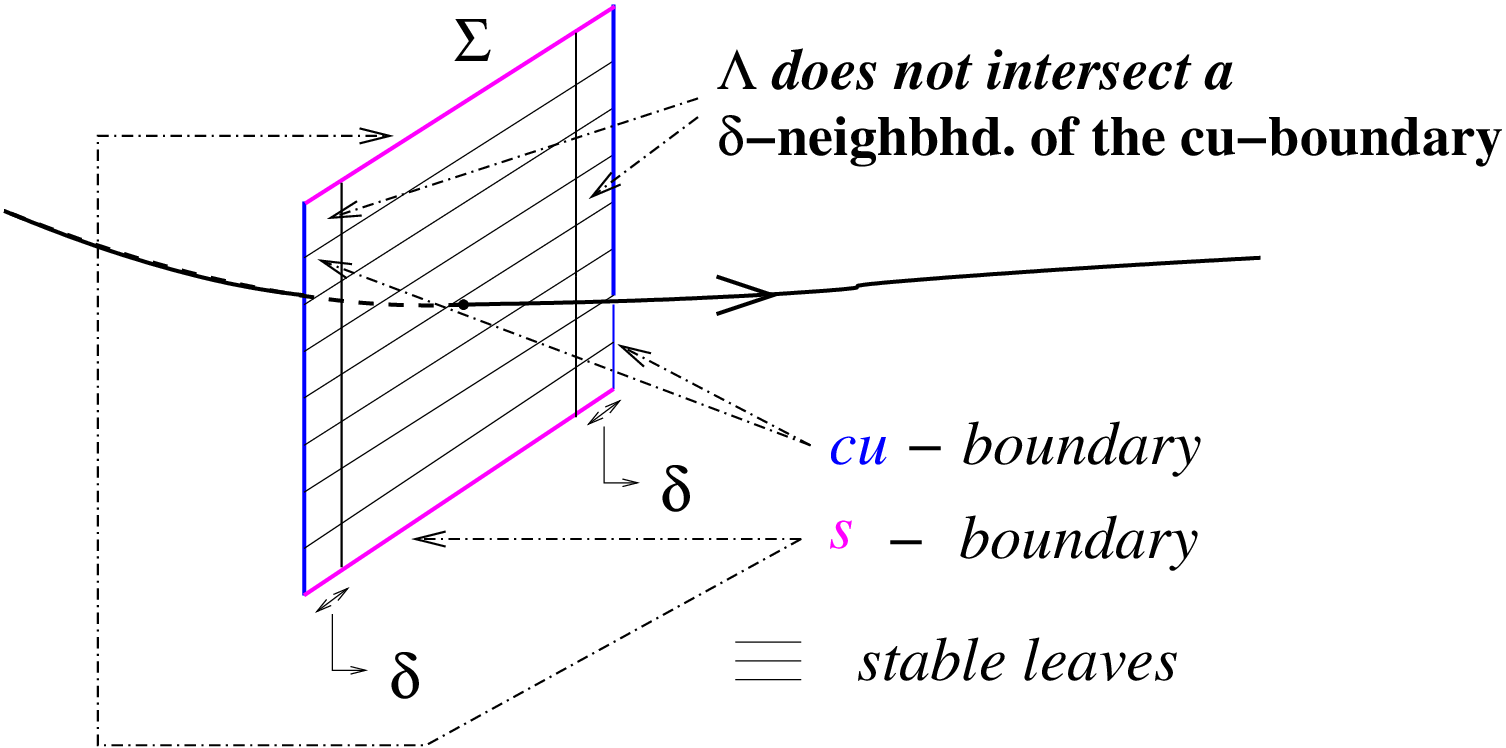}
  \caption{An adapted cross-section for $\Lambda$.}
  \label{fig:an-adapted-cross}
\end{figure}

Now we exhibit stable manifolds for Poincar\'e
transformations $R:\Sigma\to\Sigma'$. The natural candidates
are the intersections $W^s(x,\Sigma)=W^s_\epsilon(x)\cap\Sigma$
we introduced previously.  These intersections are tangent
to the corresponding sub-bundle $E^s_\Sigma$ and so, by
Proposition~\ref{p.secaohiperbolica}, they are contracted by
the transformation.  For our purposes it is also important
that the stable foliation be invariant:
\begin{equation}\label{eq.stableMarkov}
R(W^s(x,\Sigma)) \subset W^s(R(x),\Sigma')
\qquad \text{for every } x\in\Lambda\cap\Sigma.
\end{equation}
In order to have this we restrict our class of
cross-sections whose center-unstable boundary is disjoint
from $\Lambda$.
Recall  that we are considering
cross-sections $\Sigma$ that are diffeomorphic to the square
$[0,1]\times[0,1]$, with the horizontal lines
$[0,1]\times\{\eta\}$ being mapped to stable sets
$W^s(y,\Sigma)$.  The \emph{stable boundary}
$\partial^{s}\Sigma$ is the image of $[0,1]\times\{0,1\}$.
The \emph{center-unstable boundary} $\partial^{cu}\Sigma$ is
the image of $\{0,1\}\times [0,1]$. The cross-section is
\emph{$\delta$-adapted} if
$$
d(\Lambda \cap \Sigma,\partial^{cu}\Sigma)> \delta,
$$
where $d$ is the intrinsic distance in $\Sigma$, see
Figure~\ref{fig:an-adapted-cross}.

\begin{lemma}\label{l.existeadaptada}
  Let $x \in \Lambda$ be a regular point, that is, such that
  $X(x)\neq 0$. There exists $\delta>0$ such that there
  exists a $\delta$-adapted cross-section $\Sigma$ at $x$.
\end{lemma}

We are going to show that if the cross-sections are adapted,
then we have the invariance property
\eqref{eq.stableMarkov}.  Given $\Sigma,\Sigma'\in\Xi$ we
set $\Sigma(\Sigma')=\{ x\in\Sigma: R(x)\in\Sigma'\}$ the
domain of the return map from $\Sigma$ to $\Sigma'$.

\begin{lemma}
\label{stablereturnmap}
Given $\delta>0$ and $\delta$-adapted cross-sections
$\Sigma$ and $\Sigma'$, there exists
$T_2=T_2(\Sigma,\Sigma')>0$ such that if
$R:\Sigma(\Sigma')\to\Sigma'$ defined by $R(z)=R_{t(z)}(z)$
is a Poincar\'e map with time $t(\cdot)>T_2$, then
\begin{enumerate}
\item $R\big(W^s(x,\Sigma)\big)\subset W^s(R(x),\Sigma')$
      for every $x\in\Sigma(\Sigma')$, and also
\item $d(R(y),R(z))\le \frac12 \, d(y,z)$ for every $y$,
      $z\in W^s(x,\Sigma)$ and $x\in\Sigma(\Sigma')$.
\end{enumerate}
\end{lemma}
Clearly we may choose $T_2>T_1$ so that all the properties
of the Poincaré maps obtained up to here are valid for
return times greater than $T_2$.


\subsection{The proof of expansiveness}
\label{sec:proof-expansiveness}

Here we sketch the proof of
Theorem~\ref{thm:sing-hyp-attract-expansive}.  The proof is
by contradiction: let us suppose that there exist
$\epsilon>0$, a sequence $\delta_n\to 0$, a sequence of
functions $h_n\in \cK$ and sequences of points $x_n,\, y_n
\in\Lambda$ such that
\begin{equation}
\label{eq.rel1}
d\big(X^t(x_n),X^{h_n(t)}(y_n)\big)\le\delta_n
\quad\text{for all } t\in\RR,
\end{equation}
but
\begin{equation}
\label{eq.rel2}
X^{h_n(t)}(y_n) \notin X^{[t-\epsilon,t+\epsilon]}(x_n)
\quad\text{for all } t\in\RR.
\end{equation}
The main step of the proof is a reduction to a forward
expansiveness statement about Poincar\'e maps which
we state in Theorem~\ref{t.expansivepoincare} below.

We are going to use the following observation: there exists
some regular (i.e. non-equilibrium) point $z\in\Lambda$
which is accumulated by the sequence of $\omega$-limit sets
$\omega(x_n)$.  To see that this is so, start by observing
that accumulation points do exist, since $\Lambda$ is compact.
Moreover, if the $\omega$-limit sets accumulate on a
singularity then they also accumulate on at least one of the
corresponding unstable branches which, of course, consists
of regular points.  We fix such a $z$ once and for
all. Replacing our sequences by subsequences, if necessary,
we may suppose that for every $n$ there exists
$z_n\in\omega(x_n)$ such that $z_n\to z$.

Let $\Sigma$ be a $\delta$-adapted cross-section at $z$, for
some small $\delta$. Reducing $\delta$ (but keeping the same
cross-section) we may ensure that $z$ is in the interior of
the subset
$$
\Sigma_\delta=\{y\in\Sigma: d(y,\partial\Sigma)>\delta\}.
$$
By definition, $x_n$ returns infinitely often to the neighborhood
of $z_n$ which, on its turn, is close to $z$.
Thus dropping a finite number of terms in our sequences if
necessary, we have that the orbit of $x_n$ intersects
$\Sigma_\delta$ infinitely many times.
Let $t_n$ be the time corresponding to the $n$th
intersection.

Replacing $x_n$, $y_n$, $t$, and $h_n$ by
$x^{(n)}=X^{t_n}(x_n)$, $y^{(n)}=X^{h_n(t_n)}(y_n)$, $t'=t-t_n$,
and $h_n'(t')=h_n(t'+t_n)-h_n(t_n)$, we may suppose that
$x^{(n)}\in\Sigma_\delta$\,, while preserving both relations
\eqref{eq.rel1} and \eqref{eq.rel2}.  Moreover there exists
a sequence $\tau_{n,j}$\,, $j\ge 0$ with $\tau_{n,0}=0$ such
that
\begin{equation}
\label{eq.xnj}
x^{(n)}(j)=X^{\tau_{n,j}}(x^{(n)})\in\Sigma_\delta
\quad\text{and}\quad
\tau_{n,j}-\tau_{n,j-1}>\max\{t_1,t_2\}
\end{equation}
for all $j\ge 1$, where $t_1$ is given by
Proposition~\ref{p.secaohiperbolica} and $t_2$ is given by
Lemma~\ref{stablereturnmap}.

\begin{theorem}
\label{t.expansivepoincare}
Given $\epsilon_0>0$ there exists $\delta_0>0$ such that if
$x\in\Sigma_\delta$ and $y\in\Lambda$ satisfy
\begin{description}
\item[(a)] there exist $\tau_j$ such that
\[
x_j=X^{\tau_j}(x)\in\Sigma_\delta
\quad\text{and}\quad
\tau_{j}-\tau_{j-1}>\max\{t_1,t_2\}
\quad\text{for all $j\ge 1$};
\]
\item[(b)] $\dist \big(X^t(x),X^{h(t)}(y)\big) < \delta_0$, for all $t>0$
and some $h\in\cK$;
\end{description}
then there exists $s=\tau_j\in\RR$ for some $j\ge1$ such
that $X^{h(s)}(y)\in
W_{\epsilon_0}^{ss}(X^{[s-\epsilon_0,s+\epsilon_0]}(x))$.
\end{theorem}

The proof of Theorem~\ref{t.expansivepoincare} will not be
given here, and can be found in~\cite{APPV}.  We explain why
this implies Theorem~\ref{thm:sing-hyp-attract-expansive}. We
are going to use the following observation.

\begin{lemma}\label{proximo}
  There exist $\rho>0$ small and $c>0$, depending only on
  the flow, such that if $z_1, z_2, z_3$ are points in
  $\Lambda$ satisfying $z_3\in X^{[-\rho,\rho]}(z_2)$ and
  $z_2\in W_\rho^{ss}(z_1)$, with $z_1$ away from any
  equilibria of $X$, then
\[
\dist(z_1,z_3) \ge c \cdot \max\{\dist(z_1,z_2),\dist(z_2,z_3)\}.
\]
\end{lemma}

This is a direct consequence of the fact that the angle
between $E^{ss}$ and the flow direction is bounded from zero
which, on its turn, follows from the fact that the latter is
contained in the center-unstable sub-bundle $E^{cu}$.

\begin{figure}[ht]
\centering
\psfrag{a}{$X^{[-\rho,\rho]}(X^{h(s)-t}(y))$}
\psfrag{b}{$X^{h(s+\tau - t)}(y)$}
\psfrag{c}{$X^{h(s)-t}(y)$}
\psfrag{d}{$X^{h(s)}(y)$}
\psfrag{e}{$X^{s+\tau}(x)$}
\psfrag{f}{$W^{ss}_\epsilon(X^{s+\tau}(x))$}
\psfrag{g}{$X^{s+\tau-t}(x)$}
\psfrag{h}{$W^{ss}_\rho(X^{s+\tau-t}(x))$}
\includegraphics[height=3.5cm]{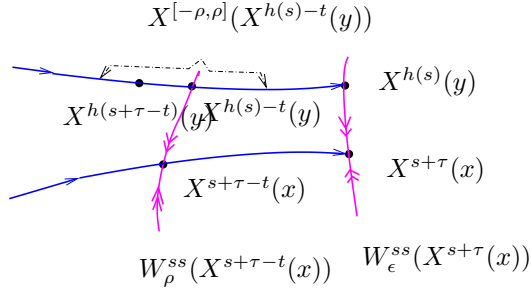}
\caption{\label{fig.contradiction} Relative
  positions of the strong-stable manifolds and orbits.
}
\end{figure}

We fix $\epsilon_0=\epsilon$ as in \eqref{eq.rel2} and then
consider $\delta_0$ as given by
Theorem~\ref{t.expansivepoincare}.  Next, we fix $n$ such
that $\delta_n<\delta_0$ and $\delta_n < c \rho$, and apply
Theorem~\ref{t.expansivepoincare} to $x=x^{(n)}$ and
$y=y^{(n)}$ and $h=h_n$.  Hypothesis (a) in the theorem
corresponds to \eqref{eq.xnj} and, with these choices,
hypothesis (b) follows from \eqref{eq.rel1}.  Therefore we
obtain that $X^{h(s)}(y)\in
W_\epsilon^{ss}(X^{[s-\epsilon,s+\epsilon]}(x))$.
Equivalently there is $|\tau|\le\epsilon$ such that
$X^{h(s)}(y)\in W_\epsilon^{ss}(X^{s+\tau}(x))$.  Condition
\eqref{eq.rel2} then implies that $X^{h(s)}(y) \neq
X^{s+\tau}(x)$.  Hence since strong-stable manifolds are
expanded under backward iteration, there exists $\theta>0$
maximum such that
$$
X^{h(s)-t}(y)\in W^{ss}_\rho(X^{s+\tau-t}(x))
\quad\text{and}\quad
X^{h(s+\tau-t)}(y)\in X^{[-\rho,\rho]}(X^{h(s)-t}(y))
$$
for all $0\le t\le \theta$, see
Figure~\ref{fig.contradiction}.  Moreover $s=\tau_j$ for
some $j\ge1$ so that $x$ is close to cross-section of the
flow which we can assume is uniformly bounded away from the
equilibria, and then we can assume that $\|X(X^t(x))\|\ge c$
for $0\le t\le \theta$.  Since $\theta$ is maximum
\begin{align*}
  \text{either  }
  \dist\big(X^{h(s)-t}(y),X^{s+\tau-t}(x)\big)&\ge\rho
  \\
  \text{  or  }
  \dist\big( X^{h(s+\tau-t)}(y), X^{h(s)-t}(y)\big)&\ge c_0\rho
\end{align*}
for $t=\theta$, because $\|X(X^t(x))\| \ge c_0>0$ for $0\le
t\le\theta$.
Using Lemma~\ref{proximo}, we conclude that
$\dist(X^{s+\tau-t}(x),X^{h(s+\tau-t)}(y)) \ge c \rho
>\delta_n $ which contradicts \eqref{eq.rel1}.  This
contradiction reduces the proof of
Theorem~\ref{thm:sing-hyp-attract-expansive} to that of
Theorem~\ref{t.expansivepoincare}.

\subsection{Singular-hyperbolicity and chaotic behavior}
\label{sec:singul-hyperb-chaoticbehav}

Here we explain why singular-hyperbolic attractors, like the
Lorenz attractor, are necessarily robustly chaotic.

\begin{proof}[of Theorem~\ref{thm:sing-hyp-chaotic}]
  The assumption of singular-hyperbolicity on an isolated
  proper subset $\Lambda$ with isolating neighborhood $U$
  ensures that the maximal invariant subsets
  $\cap_{t\in\RR}\overline{ Y^t(U)}$ for all $C^1$ nearby
  flows $Y$ are also singular-hyperbolic. Therefore to
  deduce robust chaotic behavior in this setting it is
  enough to show that a proper isolated invariant compact
  singular-hyperbolic subset is chaotic.

  Let $\Lambda$ be a singular-hyperbolic isolated proper
  subset for a $C^1$ flow. Then there exists a strong-stable
  manifold $W^{ss}(x)$ through each of its points $x$. We
  claim that this implies that $\Lambda$ is past
  chaotic. Indeed, assume by contradiction that we can find
  $y\in W^{ss}(x)$ such that $y\neq x$ and
  $\dist\big(X^{-t}(y), X^{-t}(x)\big) < \epsilon$ for every
  $t>0$, for some small $\epsilon>0$. Then, because
  $W^{ss}(x)$ is uniformly contracted by the flow in
  positive time, there exists $\lambda>0$ such that
  \begin{align*}
    \dist(y,x)\le Const\cdot e^{-\lambda
      t}\dist\big(X^{-t}(y), X^{-t}(x)\big) \le Const\cdot
    \epsilon e^{-\lambda t}
  \end{align*}
  for all $t>0$, a contradiction since $y\neq x$. Hence for
  any given small $\epsilon>0$ we can always find a point
  $y$ arbitrarily close to $x$ (it is enough to choose $y$
  is the strong-stable manifold of $x$) such that its past
  orbit separates from the orbit of $x$.

  To obtain future chaotic behavior, we argue by
  contradiction: we assume that $\Lambda$ is not future
  chaotic. Then for every $\epsilon>0$ we can find a point
  $x\in\Lambda$ and an open neighborhood $V$ of $x$ such
  that the future orbit of each $y\in V$ is $\epsilon$-close
  to the future orbit of $x$, that is, $\dist\big(X^{t}(y),
  X^{t}(x)\big) \le \epsilon$ for all $t>0$. 

  First, $x$ is not a singularity, because all the possible
  singularities inside a singular-hyperbolic set are
  hyperbolic saddles and so each singularity has a unstable
  manifold. Likewise, $x$ cannot be in the stable manifold
  of a singularity. Therefore $\omega(x)$ contains some
  regular point $z$. Let $\Sigma$ be a transversal section
  to the flow $X^t$ at $z$. 

  Hence there are infinitely many times $t_n\to+\infty$ such
  that $x_n:=X^{t_n}(x)\in\Sigma$ and $x_n\to z$ when
  $n\to+\infty$. Taking $\Sigma$ sufficiently small looking
  only to very large times, the assumption on $V$ ensures
  that each $y\in V$ admits also an infinite sequence
  $t_n(y)\xrightarrow[n\to+\infty]{}+\infty$ satisfying
  \begin{align*}
    y_n:=X^{t_n(y)}(y)\in\Sigma \quad\text{and}\quad
    \dist(y_n, x_n)\le 10\epsilon.
  \end{align*}
  We can assume that $y\in V$ does not belong to $W^s(x)$,
  since $W^s(x)$ is a $C^1$ immersed sub-manifold of
  $M$. Hence we consider the connected components
  $\gamma_n:=W^s(x_n,\Sigma)$ and $\xi_n:=W^s(y_n,\Sigma)$
  of $W^s(x)\cap\Sigma$ and $W^s(y)\cap\Sigma$,
  respectively. We recall that we can assume that every $y$
  in a small neighborhood of $\Lambda$ admits an invariant
  stable manifold because we can extend the invariant stable
  cone fields from $\Lambda$ to a small neighborhood of
  $\Lambda$. We can also extend the invariant
  center-unstable cone fields from $\Lambda$ to this same
  neighborhood, so that we can also define the notion of
  $cu$-curve in $\Sigma$ in this setting.

\begin{figure}[htpb]
\centering
\psfrag{a}{$y_j$}
\psfrag{b}{$x_j$}
\psfrag{c}{$y_{j+1}$}
\psfrag{d}{$x_{j+1}$}
\psfrag{e}{$\gamma_{j+1}$}
\psfrag{f}{$\gamma_{j}$}
\psfrag{g}{$\Sigma^{j}$}
\psfrag{h}{$\Sigma^{j+1}$}
\includegraphics[height=3cm]{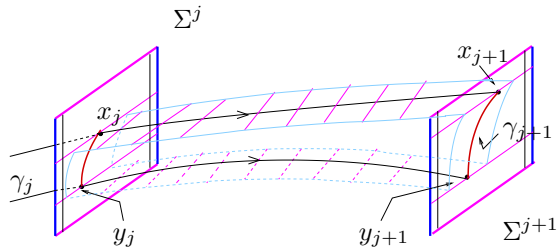}
\caption{Expansion between visits to a cross-section.}
  \label{fig:expans-within-tube}
\end{figure}

  The assumption on $V$ ensures that there exists a
  $cu$-curve $\zeta_n$ in $\Sigma$ connecting $\gamma_n$ to
  $\xi_n$, because $X^{t_n}(V)\cap\Sigma$ is an open
  neighborhood of $x_n$ containing $y_n$. But we can assume
  without loss of generality that
  $t_{n+1}-t_n>\max\{t_1,t_2\}$, forgetting some returns to
  $\Sigma$ in between if necessary and relabeling the times
  $t_n$. Thus Proposition~\ref{p.secaohiperbolica} applies
  and the Poincaré return maps associated to the returns to
  $\Sigma$ considered above are hyperbolic.

  The same argument as in the proof of expansiveness
  guarantees that there exists a flow box connecting
  $\{x_n,y_n\}$ to $\{x_{n+1},y_{n+1}\}$ and sending
  $\zeta_n$ into a $cu$-curve $R(\zeta_n)$ connecting
  $\gamma_{n+1}$ and $\xi_{n+1}$, for every $n\ge1$.

  The hyperbolicity of the Poincaré return maps ensures that
  the length of $R(\zeta_n)$ grows by a factor greater than
  one, see Figure~\ref{fig:expans-within-tube}. Therefore,
  since $y_n,x_n$ are uniformly close, this implies that the
  length of $\zeta_1$ and the distance between $\gamma_1$
  and $\xi_1$ must be zero. This contradicts the choice of
  $y\neq W^s(x)$.

  This contradiction shows that $\Lambda$ is future chaotic,
  and concludes the proof. 
\end{proof}



\def\cprime{$'$}


\end{document}